\documentclass[12pt,fleqn,a4paper]{article}
 \usepackage{amssymb}
 \usepackage{amsfonts}
 \usepackage{amsmath}
\usepackage{graphicx}
\usepackage{latexsym}
\usepackage{enumerate}
 \usepackage{dsfont}
\usepackage{color}
\usepackage{stmaryrd}
 \usepackage{mathrsfs}
 \usepackage{multicol}

 \newtheorem{theorem}{\sc Theorem}[section]
 
 \newtheorem{lemma}[theorem]{\sc Lemma}
 
 \newtheorem{corollary}[theorem]{\sc Corollary}

 \newtheorem{remark}[theorem]{\sc Remark}
 \numberwithin{equation}{section}

 \DeclareMathAlphabet{\mathsfn}{OT1}{cmss}{x}{n}
\DeclareMathAlphabet{\mathsfb}{OT1}{cmss}{bx}{n}

  \newcommand{\bfm}{\boldsymbol}
  \newcommand{\alg}{\bfm}
 
 \newcommand{\bg}{\begin}
 \newcommand{\e}{\end}
 \newcommand{\mbb}{\mathsf}
 \newcommand{\be}{\bg{enumerate}}
 \newcommand{\ee}{\e{enumerate}}
 \newcommand{\bi}{\bg{itemize}}
 \newcommand{\ei}{\e{itemize}}
 \newcommand{\bd}{\bg{description}}
 \newcommand{\ed}{\e{description}}
 \newcommand{\Ra}{\Rightarrow}
 \newcommand{\ra}{\to}
 \newcommand{\La}{\Leftarrow}

 \newcommand{\beqn}{\bg{eqnarray*}}

 \newcommand{\eeqn}{\e{eqnarray*}}
 \newcommand{\smsm}{\smallsetminus}
 
 \newcommand{\no}{\noindent}
 
 \newcommand{\A}{\alg{A}}
 \newcommand{\B}{\alg{B}}

 \newfont{\gros}{cmr10 at 13pt}

\def\le{\leqslant}
\def\ge{\geqslant}

\def\brl{boun\-ded resi\-duated lattice}
\def\imf{implicative filter}
\def\brls{boun\-ded re\-si\-duated lattices}
\def\V{\mbb{V}}
\def\Vk{\mbb{V^k}}
\def\F{\bfm{F}}

\def\BRL{\mbb{BRL}}
\def\WL{\mbb{WL}}
\def\WLk{\mbb{WL_\mbb{k}}}

\def\SI{\mbb{SI}}

\def\Em{\mbb{E_{m}}}
\def\EMn{\mbb{EM}_{\mbb{n}}}
\def\EMm{\mbb{EM}_{\mbb{m}}}

\def\cqd{\hfill$\Box$}

 \newcommand{\LLn}{L^{\omega}_{n+1}}

\begin{document}

 \title{The General Apple Property and Boolean terms in Integral Bounded Residuated Lattice-ordered Commutative Monoids.}
\author{by \sc Antoni Torrens}
\date{}
\maketitle

\begin{abstract}

 In this paper we give  equational presentations of the varieties of {\em integral bounded residuated lattice-ordered commutative monoids}
  (bounded residuated lattices for short) satisfying the \emph{General Apple Property} (GAP),  that is,  varieties in which  all of its directly indecomposable members are local. This characterization is given by means of Boolean terms: \emph{A variety $\V$ of \brl s has GAP iff there is an unary term $b(x)$  such that $\V$ satisfies the equations $b(x)\lor\neg b(x)\approx \top$ and $(x^k\to b(x))\cdot(b(x)\to k.x)\approx \top$, for some $k>0$}. Using this characterization, we show that  for   any variety $\V$ of bounded residuated lattice satisfying GAP there is $k>0$ such that the equation $k.x\lor k.\neg x\approx \top$ holds in $\V$, that is, $\V\subseteq \WLk$. As a consequence  we improve Theorem 5.7 of  \cite{CT12}, showing in theorem that a\emph{ variety of \brls\ has Boolean retraction term if and only if there is $k>0$ such that it satisfies the equation  $k.x^k\lor k.(\neg x)^k\approx\top$.}
  We also see that in Bounded residuated lattices GAP is equivalent to Boolean lifting property (BLP) and so, it is  equivalent to quasi-local property (in the sense of \cite{GLM12}). Finally, we prove that a variety of \brl s has   GAP and its semisimple members form a variety if and only if there exists an unary term which is simultaneously Boolean and radical  for this variety.
 \end{abstract}


 \section*{Introduction}

  The Apple property (AP) was considered in \cite{BB87} to study  free algebras in  locally finite varieties having Fraser-Horn property (FHP).
  Following \cite{BB87} a finite algebra $\A$ has the \emph{Apple Property} (AP) if for any proper factor congruence relation $\theta$ on $\A$, $\theta\not=A\times A$,  such that
 $\A/\theta$ is directly indecomposable,  then the interval lattice $[\theta, A\times A]$ has exactly one coatom of $\langle Con(\A),\subseteq\rangle$, that is, the algebra $\A/\theta$ is \emph{local}, i.e., it has a unique maximal proper congruence relation.
 A variety $\V$ has AP if  any $\A\in \V$  finite and directly indecomposable is local,
 or equiva\-lently, if any finite algebra in $\V$ has AP.\smallskip

 In this paper we consider a more general condition.
 We say that an algebra $\A$ has the \emph{General Apple Property} (GAP) if for
  any proper  congruence relation $\theta$ on $\A$ such that
 $\A/\theta$ is directly indecomposable, then  $\A/\theta$ is local, and a variety $\V$ has GAP when any algebra in $\V$ has GAP; i.e, (see Lemma \ref{L:GAP-Locrep}), $\V$ has GAP if and only if  any directly indecomposable in $\V$  is local.
 \smallskip

 Our study is focused in integral bounded residua\-ted lattice-ordered commutative monoids (bounded residuated lattices for short). They form an  arithmetical variety $\BRL$, i.e., its members are congruence distributive and congruence permutable, and so they have FHP.\smallskip

 For a \brl\ $\A$, GAP can be stated in the next form:
 \begin{quote}\em
 ``If $F$ is a proper implicative filter such that $\A/F$ is directly indecomposable, then  $F$ is contained in a unique maximal proper implicative filter".
 \end{quote}

 By definition, a {\em factor congruence relation}  on an algebra  is a congruence relation on this algebra complemented in the congruence lattice, and which permutes with one of its complements. Algebras whose factor congruences  form a Boolean algebra are said to have the {\em Boolean factor congruence pro\-perty} (BFP). Varieties with FHP have BFP (see \cite{BiBu90}), and so, varieties of \brl s have BFP.
  Since in bounded  residuated lattices factor congruences correspond with Boolean (complemented) elements, we deduce (see Lemma~\ref{GAP-Stone}) that in any \brl\ $\A$ GAP is equivalent to:
 \begin{quote}\em
 ``Every Stone ultrafilter of $\A$ is contained in a unique maximal implicative filter of $\A$",
 \end{quote}
 where by Stone ultrafilter we mean  an implicative filter generated by an ultrafilter of the Boolean algebra of Boolean elements of $\A$.\smallskip

 Using this fact, we give in Theorem \ref{T:rad(f)-max} a characterization of GAP for \brl s in si\-mi\-lar way as it is done in \cite[Theorem 3.3]{CT96} for MV-algebras. From this result we deduce that GAP is equivalent to Boolean Lifting Property (BLP), in the sense of \cite{GM14}.
 In Theorem \ref{T:GAP_bt} we give equational presentations of the varieties of bounded residuated lattices satisfying GAP, by means of \emph{Boolean terms} defined in Section 2. As a consequence, we obtain that  each variety of \brl s satisfying GAP  is contained in $\WLk$ for some   $k>0$, that is,  it satisfies the equation $k.x\lor k.\neg x\approx\top$. We also characterize the  varieties of \brl s satisfying GAP  whose semisimple members form a variety, by means of the existence of  radical and Boolean terms. As a consequence, in Theorem
\ref{T:Br-WLn}, we improve Theorem 5.7. of  \cite{CT12}, showing that a variety of \brls\ has\emph{ Boolean retraction term} if and only if there is $n>0$ such that it satisfies the equation  $n.x^n\lor n.(\neg x)^n\approx\top$.   These results are described in Section 3, which is the main section of the paper.
\smallskip

  To obtain these results it was necessary to study the Boolean algebra formed by complemented elements in \brl s, the implicative filters gene\-rated by its Boolean filters (ultrafilters), called Stone filters (Stone ultrafilters), and its relation with  direct decomposability. This is done in Section 2, in which we  also introduce  Boolean terms,  needed  in section 3, giving some properties of them.\smallskip

 The paper also contains a preliminary  section in which we give the notation and some results on bounded  residuated lattices needed in the paper. We explain them in some detail so that the paper is as self-contained as possible; however, more information on these results can be found in \cite{GJKO,CT12,T16} and in the papers and books referenced therein.
 This section has two subsection, namely  1.1 and 1.2.
 In subsection 1.1, we recall some  arithmetical properties, the family of varieties  $(\WLk)_{k>0}$  introduced in \cite{CT12} and the family  $(\Em)_{m>0}$  of $m$-contrac\-ti\-ve varieties $m>0$, which play an important role in the paper.
  In subsection 1.2, we  recall the relation between implicative filters and congruence relations, the most common properties on prime filters, maximal filters,  and their relation with directly indecomposable, simple and semisimple members of $\BRL$.
  We also recall some properties on semisimplicity and its relation with radical term.\smallskip

  Finally, at the end of the paper, we present two examples of varieties of \brl s admitting boolean term satisfying the condition (\ref{Eq:VDI=VL1}) and so, by Theorem \ref{T:GAP_bt},  they are GAP. Furthermore, both admit a Boolean term, which is also a radical term, and thus the classes of their semisimple members are varieties.\smallskip

 The reader familiar with \brl s  only need a quick reading of section 1, paying special attention to the nomenclature and notation used in the paper. We  assume that the reader is also  familiar with the  basic results on Universal algebra which can be found in the books \cite{BS81}, \cite{Be12} and in the references given therein.

 \section{Preliminaries}
 \subsection{Arithmetical properties of \brls.}
 Throughout this paper $\mbb{BRL}$ will denote the class of all {\em bounded integral residuated lattice-ordered commutative  monoid (bounded residuated lattices for short)}, that is, the class of algebras
 $\A  =\langle A; \cdot , \to, \lor, \land, \top, \bot\rangle $
 in the algebraic language $\{\cdot, \to,\land, \lor,  0, 1\}$ of type $(2,2,2,2,0,0)$
  such that:
   \bi
   \item $\langle A; \cdot , \top\rangle$ is a commutative monoid,
 \item $\langle A; \lor, \land,\bot, \top\rangle$ is a bounded
 lattice with smallest element $\bot$, and greatest element $\top,$
 \item  It satisfies the following residuation
 condition,  for any $a,b,c\in A$
 \bg{align}\label{Eq:resd-cond}
 a \cdot  b \le  c &
 \mbox{ if and only if } a \le b\to c,\end{align}
 where  $\le $ is the order given by the
 lattice structure.
 \ei

  It is well known that \brls\ admit  an equational presentation, and so $\BRL$ is a variety (see for example \cite{KO01,GJKO} and \cite{CT12}), that is,  if
 $H$, $S$ and $P$ represents respectively the  class operator homomorphic images, isomorphic images of subalgebras and isomorphic images of direct products, then $HSP(\BRL)\subseteq \BRL$. In fact $\BRL$ is an arithmetical variety, i.e., congruence permutable and congruence distributive.
  \smallskip

 For further references we list   some well known pro\-perties of \brl s. Let $\A\in \BRL$, then
  \bg{lemma}\label{basic}
  The following properties hold true
  for any $a,b,c\in A$:

 \bg{enumerate}[$a)$]\begin{multicols}{2}
  \item\label{order} $a \le  b $
  if and only if $a \to b = \top$,
  \item\label{timx} $\top \to a =
  a$,
  \item\label{intercambio} $(a \cdot  b) \to c= a \to (b \to
  c)$,
  \item \label{DeMorgan} $(a \lor b) \to c=(a \to c) \land (b \to c))$,
  \item \label{resi-and} $a\to (b\land c)= (a\to b)\land (a\to c)$,
  \item \label{dis do-or} $a\cdot(b\lor c)=(a\cdot b)\lor(a\cdot c)$.
\end{multicols}
\e{enumerate}
  \e{lemma}

 We consider the  unary operation term
 $\neg x =: x \to
 \bot$, then, 
 by taking into account that the $\{\to,\bot,\top\}$-reduct
 of a bounded residuated lattice is a bounded BCK-algebra, we have (\cite{GT07,CT12}):
 \bg{lemma}\label{basicbounded}  The following properties hold true for any $a,b\in A$:
\bg{enumerate}[$a)$]\begin{multicols}{2}
 \item\label{bb1} $a \le  b\Ra\neg b \le  \neg a$,
 \item\label{bb3} $a \le \neg \neg a   $,
 \item\label{bb2} $\neg a = \neg \neg \neg a$,
 \item\label{bb4} $a \ra \neg b = b \ra \neg a$,
 \item\label{bb5} $a \ra \neg b = \neg \neg a \ra \neg b$,
 \item\label{bb6} $\neg \neg(a \ra  \neg b) = a \ra  \neg b$,
 \item\label{bb7} $\lnot a\to (a\to b)=   \top $.
 \item \label{bb8} $a\cdot b=\top$ iff $a=\top$ and $b=\top$\hfill $\Box$
\end{multicols}\end{enumerate}
 \end{lemma}
 If we consider the binary term operation
 $x +  y =:\neg(\neg x\cdot  \neg y)$,
 then $\langle A; + \rangle$ is a commutative semigroup. We also consider the terms:
 \bi
 \item $x^0=\top$ and $0.x=\bot$,
 \item $x^{n+1}=x\cdot  x^n$ and $(n+1).x=x +  n.x$ for all nonnegative integer $n$.\ei
 Then it is straightforward to prove the  next lemma (see \cite{CT12} for example).
 \bg{lemma}\label{L:nmplus}
 For any $a,b\in A$ and for any  $0<n,m\in\omega$ it satisfies:
 \begin{multicols}{2}\setlength{\columnsep}{-2.1in}\be[$a)$]
  \item\label{plus0} $a +  b=\neg a\to \neg\neg b$,
  \item\label{plus1} $1.a=\bot +a=\neg\neg a$,
  \item\label{plus2} $n.a =\neg(\neg a)^n$,
  \item\label{plus3} $n.a = \neg\neg (n.a)=n.(\neg\neg a)$,
  \item\label{plus4} $(n+m).a= (n.a) +  (m.a)$,
  \item\label{plus4b} $n.(a+b)=n.a+n.b$,
  \item\label{plus4a} $(m\cdot n).a = m.(n.a)$,
  \item\label{plus5} if $n\le m$, then \;$n.a\le m.a$ and $a^m \le  a^n$,
  \item\label{plus6}  $\neg((n.a)^m)=m.(\neg a)^n$,
  \item\label{plus61} $n.(\neg(x^m))= \neg (x^{mn})$,
  \item\label{plus7}  if $a\le b$, then $n.a^m\le n.b^m$.\cqd
 \ee \end{multicols}
 \e{lemma}

 Let $n>0$. An element $a\in A$ is called $n$-\emph{contractive} or  $n$-\emph{potent} when $a^n=a^{n+1}$; $2$-contractive elements are also known as  \emph{idempotent}.
 Then $\A$ is called \emph{locally contractive} if any element  $a\in A$  is $m$-contractive for some $0<m\in\omega$.
\smallskip

 For any  integer $m>0$, let $\Em$ denote the subvariety of $\BRL$  determined by  the identity
   \bd
     \item[ \rm  \ (Em)] \ $x^m\approx x^{m+1}$.
   \ed
   That is, the members of $\Em$ are those that all its elements are $m$-contractive.
  It is well known that for all $m>0$, $\mbb{E_m}\subsetneq\mbb{E_{m+1}}$. Algebras in  $\mbb{E_1}$ are also called  \emph{Heyting algebras}\label{D:Heyting}.

  \bg{lemma} Let $\V$ be a variety of \brl. Then the following conditions are equivalent:
   \be[$(i)$]
   \item Any member of $\V$ is locally contractive.
   \item The 1-free algebra of $\V$ is locally contractive.
   \item $\V\subseteq\Em$, for some $m>0$.
   \ee
  \e{lemma}
  {\em Proof:}
  It is plain that  $(i)$ implies $(ii)$ and $(iii)$ implies $(i)$. $(ii)$ implies $(iii)$, because $(ii)$ implies that the free generator of the $1$-free algebra of $\V$ is $m$-contractive for some $m>0$, and so $\V\models x^m\approx  x^{m+1}$.\hfill$\Box$

  An algebra is called \emph{locally finite} when  its finitely generated subalgebras are finite. Varieties whose members are locally finite are called \emph{locally finite varieties}. Since  $1$-free algebras in  locally finite varieties   are  finite, we have:
  \bg{corollary}\label{C:lf->m-cont} Any locally finite variety of \brls\ is contained in $\Em$ for some $m>0$.\hfill $\Box$
  \e{corollary}

   For each integer $k>0$,  $\mbb{WL_k}$ represents the subvariety of $\BRL$ given by the identity (see \cite{CT12})
 \bg{description}
 \item[ \rm \ (WLk)] \ $k.x\lor k.\neg x\approx \top$.\label{E:WLk}
 \e{description}
   The variety $\mbb{WL}_2$ contains the class $\mbb{MTL}$ of  MTL-\emph{algebras}, or \brls\ representable as subdirect product totally ordered bounded residuated lattices.  In fact $\mbb{MTL}$  is the  subvariety of $\BRL$ given by the identity   $(x\to y)\lor (y\to x)\approx\top$, which contains de class $\mbb{MV}$ of all MV-algebras.\smallskip

  From the results given in \cite{CT12} we deduce:
 \bg{lemma} The following properties hold  true:
 \be[$1)$]
 \item $\mbb{WL}_1$ is the variety of Stonean residuated lattices, i.e., the subvariety of $\BRL$ given by the equation $\neg\neg x\lor \neg x\approx\top$.
 \item For any integer $k>0$, $\mbb{WL}_k\varsubsetneq \mbb{WL}_{k+1}$.
 \item For any integer $k>0$, $\mbb{WL}_k\cap \mbb{PRL}=\mbb{WL}_1$, where $\mbb{PRL}$ denotes the varie\-ty of pseudocomplemented residuated lattices, i.e., the subvariety of $\BRL$ given by the identity $x\land\neg x\approx\bot$.\hfill $\Box$
 \ee
 \e{lemma}

  In what follows we  consider unary  $\{\cdot , \to,\land, \lor, \bot, \top\}$-terms called  \emph{unary terms} for short. Given an unary term $t$  we write $t(x)$ to indicate that the variable which appears in $t$  is $x$. If $\A$ is a \brl, then for any $a\in A$, $t^{\A}(a)$ represents the interpretation of $t$ on $\A$ given by the assignment $x \mapsto a$.

  Given $t,t'$ unary terms,  for notational convenience, we will  write $t\preccurlyeq t'$ in place of the equation $t\to t'\approx \top$.

  \subsection{Implicative filters, congruence relations and radical. Semisimplicity.}
 The variety $\BRL$ is $\top$-regular, that is, if $\A\in \BRL$, then every congruence relation on $\A$ is characterized by the equivalence classe of $\top$. The equivalence classes of $\top$ are called implicative filters.\smallskip

 In fact, an {\em implicative filter\/}   of a \brl\ $\A$ may be defined as a subset $F $ of $A$ satisfying the
 following  conditions:
 \be[f1)]
 \item $\top \in F$,
 \item for all $a, b\in A$, if $a \in F$ and $a \le  b$,   then $b \in
 F$, and
 \item  if $a, b$ are in $F$, then $a \cdot  b \in F$.\ee
 Alternatively, \imf\  can also  be defined as a subset $F$ of $A$
 satisfying f1) and f4) for all $a,b\in A$  $a,a \to b\in F$ imply $b \in F$.
\smallskip

 It is well known, and easy to prove, that for each non empty set $B\subseteq A$,
  \bg{equation}\label{Gen.Fil}
  \langle B\rangle=\{a\in A: x_1^{n_1}\cdot x^{n_2}\cdots x_k^{n_k}\le a,\ 1\le k,n_1,\ldots,n_k,\  x_1 ,\ldots, x_k\in B \}
  \e{equation}
  is  the \emph{smallest \imf\ of $\A$ containing} $B$, which is just the  intersection of
  all \imf s  containing $B$. For each $a \in A$, we shall write $\langle a\rangle$ instead of  $\langle \{a\}\rangle$. The next lemma is a consequence of (\ref{Gen.Fil}).
  \bg{lemma} If $\A$ is a \brl\ and  $a\in A$, then:
  \be[$(a)$]
  \item $\{\top\}$ is the least implicative filter, and $\langle \bot \rangle= A$.
  \item For any $B\subseteq A$, $\langle B\cup\{a\}\rangle=\{b\in A: a^n\to b\in \langle B\rangle, \mbox{ for some }n\ge 0\}$
  \item $\langle a\rangle=\{b\in A: a^n\to b=\top \mbox{ for some }n\ge 0\}$.\cqd
  \ee
  \e{lemma}

  Given an \imf\ $F$ of a   $\A\in \BRL$,  the binary
  relation
  \bg{equation}\label{D:ConF}
  \theta(F) := \{(x,y) \in A \times A : x \to y \in F\, \mbox{
  and  }\, y\to x \in F\}
  \e{equation}
  is a congruence on $\A$ such that $F=\top/\theta(F)$,
  the equivalence class of $\top.$
  Actually, the  correspondence  $F \mapsto \theta (F)$ is an order isomorphism from the  set of implicative filters of $\A$ onto the set of congruences of  $\A$,  both  ordered by inclusion, whose inverse is given by the map $\theta\mapsto \top/\theta$.
   We  write  $\A/F$
  instead of $\A/\theta(F)$, and $a/F$  instead of $a/\theta(F)$,   the equivalence
  class of $a$ modulo $\theta(F)$. Notice that in $\A$, $\theta(\{\top\})=\Delta_A$ is the identity relation on $A$ and $\theta(\langle\bot\rangle)=A\times A=\nabla_{A}$ is the universal equivalence relation on $A$.  \smallskip

    Implicative filters generated by one element are called  \emph{principal}.  Clearly, in a finite \brl\ any implicative filter is principal.

 An \imf\ $F$ of a non-trivial \brl\ $\A $ is called {\em proper\/} when $F \neq A$, that is,  $\bot\notin A$. A proper   implicative filter $P$ of  $\A$ is called {\em prime} when  for any $a,b\in P$, $a\lor b\in P$ implies $a\in P$ or $b\in P$. It is well known that a \emph{proper implicative filter $F$ of $\A$ is prime if and only if $\A/F$ is a non-trivial finitely subdirectly irreducible}.\label{P:FSI-Prime} If $Sp(\A)$ denotes the set of all prime \imf\ of $\A$, then, since  $\bigcap\limits_{P\in Sp(\A)}  \theta(P)=\Delta_{\A}$, the correspondence $a\mapsto \left(a/P\right)_{P\in Sp(\A)}$ gives a representation of $\A$ as  subdirect product of the family $(\A/P)_{P\in Sp(\A)}$.
 \smallskip

 Given $\A\in\BRL$, if $Spm(\A)$ denotes the set of minimal elements in the poset $\langle Sp(\A),\subseteq\rangle$, then, since  $Sp(\A)$ is closed under intersection of descen\-ding chains, we have 
 \bg{lemma}\label{L:Spm}
 Every $\A\in \BRL$ satisfies
 \be[$(a)$]
 \item For any $P\in Sp(\A)$ there is $p\in Spm(\A)$ such that $p\subseteq P$.
 \item $\bigcap Spm(\A)=\{\top\}$
 \item $\A$ is ismophic to a subdirect product of the family $(\A/p)_{ p\in Spm(\A)}$. 
     \ee
 \e{lemma}
 Given a class $\mbb{K}$ of algebras, we represent by $\mbb{K_{FSI}}$
 the class of its finitely subdirectly irreducible members and by  $\mbb{K_{SI}}$   the class of its subdirectly irreducible members. Clearly $\mbb{K_{\SI}}\subseteq\mbb{K_{FSI}}$.   Every
 variety is generated by its (finitely) subdirectly irreducible members. After \cite[Proposition 1.4]{KO01}, we know that  \emph{ $\mbb{BRL_{\mbb{FSI}}}$  is just the class of \brl s having $\top$   join irreducible.}\smallskip

  A {\em maximal \imf} of $\A$ is a proper \imf\ $M$  such that for each
  \imf\ $G$ of $\A $, $M\varsubsetneq G$ implies $G=A$. Clearly, any maximal \imf\ is prime. We denote by $Max(\A)$ the set of all maximal \imf\ of $\A$, then $Max(\A)\subseteq Sp(\A)$.
  \bg{remark}\label{R:Proper}
  Since in a \brl\ the set of its proper
   \imf s is closed under upward directed families,
 by Zorn's Lemma, an \imf\ is proper if and only if   it is contained in a maximal \imf.
 \e{remark}

The \emph{radical} of a \brl\  $\A$, represented by $Rad(\A)$, is the intersection of its maximal \imf s, that is,
 \bi
 \item $a\in Rad(\A)$ iff $a\in M$ for each maximal \imf\ $M$ of $\A$.
 \ei
 The following two results are well known and they can be found in the lite\-rature (see for example \cite{T16} and the references given therein).

 \bg{lemma}\label{L:fMax}
 Let $F$ be an \imf\ of a \brl\ $\A$. Then $F$  is maximal if and only if
 \medskip

\no \emph{\rm \ (Mx)}   for any $a\in A$, $a\notin F$ if and only if there is $n>0$  such that \ $\neg a^n\in F$.

 \e{lemma}

\bg{corollary}\label{C:nn->Max} Let $M\in Max(\A)$. Then for any $a\in A$, $a\in M$ if and only if $\neg\neg a\in M$.
\e{corollary}
{\em Proof:} If $a\in M$ then, since  $a\le\neg\neg a$, $\neg\neg\,a\in M$.
If $a\notin M$, then there is $n>0$ such that  $\neg(\neg\neg a)^n=\neg a^n\in M$, and so $\neg\neg a\notin M$.\hfill $\Box$
  \bg{lemma}[{cf. \cite[Lemma 2.4]{T16}}]\label{L:Rad}
    For every \brl\ $\A$ and for any $a\in A$, the following are equivalent.
 \be[$(i)$]
 \item $a\in Rad(\A)$.
 \item  for all $n>0$ $\langle\neg a^n\rangle =A$.
 \item for any $n>0$ there is $k>0$ such that $k.a^n=\top$\ee
 \e{lemma}

 \bg{corollary}\label{C:Def-rad} For each \brl\ $\A$ we have
 \bi
 \item[] $Rad(\A)=\{a\in A: \forall n\ge 0,\, \exists k>0 \mbox{ such that } k.a^n=\top\}$.\hfill $\Box$\ei
 \e{corollary}
 
As  a consequence of the above we have  (see \cite[Lemma 2.6]{T16}).
 \bg{lemma}\label{L: SbProdH}
  The following properties hold true:
 \be[$(a)$]
 \item If $\B$ is a subalgebra of a  $\A\in \BRL$,  then $Rad(\B)\subseteq Rad(\A)$
 \item If $(\A_i)_{i\in I}$ is a family of $\BRL$, then   $Rad(\prod\limits_{i\in I} \A_i)\subseteq \prod\limits_{i\in I} Rad( \A_i)$.
 \item If $h\colon \A\to \B$ is an homomorphism of \brls, then $h[Rad(\A)]\subseteq Rad(\B)$.
 \item For any implicative filter $F$ of   $\A\in\BRL$, $Rad(\A)/F\subseteq Rad(\A/F)$
 \ee
 \e{lemma}

  Given an integer   $k>0$, we say that an algebra  $\A\in\BRL$ is \emph{$k$-radical} provided that \begin{equation}\label{Eq:Def,k-rad}
     Rad(\A)=\{a\in A: \forall n>0,\, k.a^n=1\}.
 \end{equation}
 A variety  $\mbb{V}$  is a called \emph{$k$-radical} whenever  all its members are  $k$-radical.
 It is a simple verification that for any integer $n>0$, $\mbb{E_m}$ is $m$-radical, because for $\A\in \mbb{E_m}$, $Rad(\A)=\{a\in A: m.a^m=\top\}$. It is shown in \cite[Lemma 1.8]{CT12} that for each $k>0$,  $\mbb{WL_k}$ is    $k$-radical variety. So $\mbb{MTL}$  is $2$-radical variety, because it is a subvariety of  $\mbb{WL_2}$. \smallskip

 For $k$-radical varieties we can improve item (b) of Lemma \ref{L: SbProdH}:
 \bg{lemma}[{\cite[Lemma 2.8]{T16}}]\label{L:Prod-k-rad} Let  $k>0$.  If
$(\A_i)_{i\in I}$ is a family  of $k$-radical \brls, then $Rad(\prod\limits_{i\in I} \A_i)= \prod\limits_{i\in I} Rad( \A_i)$.\cqd
 \e{lemma}
 
An algebra  is called \emph{simple} provided that it only has trivial congruence relations, namely the identity and the universal. Hence a \brl\ $\A$ is simple if and only if $\{\top\}$ is its unique proper \imf, or equivalently, $\{\top\}$ is  maximal implicative filter. Therefore a proper \imf\  $F$ of  \emph{$\A$ is maximal if and only if $\A/F$ is simple}.
Algebras representable as subdirect product of simple algebras are called \emph{semisimple}. Therefore an  algebra is semisimple  if and only if the intersection of its proper maximal congruences  is the identity, and so  a \brl\ \emph{$\A$ is semisimple if and only if $Rad(\A)=\{\top\}$}. Notice that for any  $\A\in\BRL$ the quotient algebra $\A/Rad(\A)$ is always semisimple.
 \smallskip

 Given a variety $\V$, the class of its simple members is denoted by $\mbb{V_S}$ and  the class of its semisimple members is denoted by $\mbb{V_{SS}}$.
 It follows from Lemma \ref{L: SbProdH} that bounded residuated lattices are
 hereditarily semisimple, then 
 $\mbb{V_{SS}}$ is closed under  isomorphic images, subalgebras and products. Moreover, since $\mbb{V_{SS}}\subseteq SP(\mbb{V_S})$ and $\mbb{V_S}\subseteq \mbb{V_{SS}}$,  we have that    $\mbb{V_S}$ and $\mbb{V_{SS}}$ generate the same variety $HSP(\mbb{V_S})$.
 Notice that ${HSP(\mbb{V_S})}_{\mbb{S}}=\mbb{V_S}$, and ${HSP(\mbb{V_S})}_{\mbb{SS}}=\mbb{V_{SS}}$. \smallskip

 For any $m>1$ it is straightforward to see that $\mbb{E_{m SS}}$ is the variety $\EMm$   of \brls\ given by the equation:
 \bd
 \item[\ (EMn)] \qquad  $x\lor \neg\, x^{m}\approx \top$
 \ed
 Notice that $\mbb{EM_1}$ \emph{is the variety of Boolean algebras}.\smallskip

 We say that a variety  \emph{$\V$ is semisimple} provided that all its members are semisimple, that is, $\V=\mbb{V_{SS}}$.
 The next result follows from those  given in  \cite{K,GJKO}
    (cf.  \cite{T10}).
  \bg{theorem}\label{T:PrevioSS} Let $\mbb{V}$ be variety of bounded residuated lattices, then
   \be[\em (I)]
   \item  If $\V$ is a semisimple, then  $\V\subseteq \mbb{E_m}$ for some $m>0$,
   \item   $\V$ is  semisimple if and only if $\V\subseteq \mbb{EM_m}$, for some $m>0$,
  \item $\V$ is  semisimple  if and only if $\V$ is  discriminator variety.\hfill $\Box$\ee
  \e{theorem}

 Let $t(x)$ be  unary term. We say that a \brl\  $\A$ has $t(x)$  as \emph{radical term}, or  $t(x)$ \emph{is a radical term for} $\A$, whenever
 \bg{equation} \label{E:TermRad}
 Rad(\A)=\{a\in A: t^{\A}(a)=\top\}.
 \e{equation}
  A variety  $\mbb{V}\subseteq \BRL$  has  $t(x)$ as \emph{radical term}, or  $t(x)$ \emph{is a radical term for} $\V$, whenever $t(x)$ is radical term for each $\A\in \mbb{V}$. From Lemma\ref{L:Prod-k-rad} we deduce:

\bg{lemma}[{\cite[Lemma 5.1]{T16}}]\label{L: SbProdH2}If
$(\A_i)_{i\in I}$ is a family  of \brl s having $t(x)$ as radical term, then $Rad(\prod\limits_{i\in I} \A_i)= \prod\limits_{i\in I} Rad( \A_i)$ and $t(x)$ is also radical term for $\prod\limits_{i\in I} \A_i$.\cqd
 \e{lemma}

 Observe  that for any $m>0$ $m.x^m$ is a radical term for $\Em$, furthermore $\mbb{E_{mSS}}= \EMm$.\smallskip

  We also know that any subvariety of $\BRL$ having $t(x)$ as  Boolean retraction term (see page \pageref{P:BRT} for definition)  has   $t(x)$ as radical term. In particular, since for any $k>0$ the variety $\mbb{V^k}$ given by the equation $k.x^k\lor k.(\neg x)^k\approx \top$ has $k.x^k$ as Boolean retraction term, and then it has $k.x^k$ as radical term, see  \cite{CT12} for details. Actually, the variety $\mbb{V^k}$ is the greatest subvariety of $\WLk$ admitting Boolean retraction term, see again \cite{CT12}. In Theorem \ref{T:Br-WLn},  we will prove that any variety admitting Boolean retraction term is a subvariety of $\WLk$ for some $k>0$. Moreover, if the variety  $\V$ has  Boolean retraction term, then  $\mbb{V_{SS}}$  is the variety of Boolean algebras $\mbb{EM_1}$. 
  \smallskip

 The following theorem follows from the results given in \cite[Section 5]{T16}.

\bg{theorem}\label{T:Qe-k-rad} Let $\V$ be  a $k$-radical variety of \brls. Then
 $\V$ has  radical term if and only if $\V_{\mbb{SS}}$ is a quasivariety. In this case a radical term is  $k.x^r$, for some $r\ge k$.\hfill $\Box$
 \e{theorem}
 \bg{corollary}\label{C:wlk-trad} Let $\V$ be   a subvariety of $\mbb{WL}_k$, then
 $\V$ has  radical term if and only if $\V_{\mbb{SS}}$ is a variety. In this case a radical term is  $k.x^r$, for some $r\ge k$. \hfill $\Box$
 \e{corollary}
 \bg{remark}\label{R:krads>r} It follows from the above  that if $k.x^r$ is a radical term  of a  variety $\V$, then for any $s\ge r$ $k.x^s$ is also radical term for $\V$.
 \e{remark}

\section{Boolean elements, Stone filters, Boolean terms and indecomposability.}
In general an algebra $\A$ is called \emph{directly indecomposable} if $A$ has more than  one element and whenever it is isomorphic to a direct product of two algebras $\A_1$ and $\A_2$, then either $\A_1$ or $\A_2$ is the trivial algebra with just one element. Equivalently, see \cite[Corollary II,7.7]{BB87} the only factor congruences on $\A$ are $\Delta_A$ and $\nabla_A$.
Given a class $\mbb{K}$ of algebras we represent by $\mbb{K_{DI}}$ the class of all directly irreducible members of  $\mbb{K}$.\smallskip

 Let $\A$ be a \brl. An element $b$ of $A$ is called \emph{Boolean} if it is complemented in $\langle A,\land,\lor,\bot,\top\rangle$, in this case the complement of $b$ is $\neg b$.  The set $B(\A)$ of all complemented elements of $\A$ is universe of a subalgebra $\B(\A)$ of $\A$ which  is a Boolean algebra satisfying  $x\land y\approx x\cdot y$ and $x\lor y=x+ y$. Observe that $b\in B(\A)$ if and only if $b\lor \neg b=\top$.
 The  properties listed in following lemma are needed for further results.
\bg{lemma}\label{L:previUltra} Let $\A$ be a \brl. Then for any $a,c\in A$ and any $b\in B(\A)$, it satisfies:
\be[(a)]
\item $\neg\neg b=b$, $b\cdot a=b\land a$, $b\to a=\neg b\lor a$ and $b\land(a\lor c)=(b\land a)\lor (b\land c)$.
\item $ (b\land a)\lor \neg b=\neg b\lor a$.
\item For any $m,n>0$ $m.b^n=b$.
\ee
\e{lemma}
{\em Proof:} the properties of $(a)$ are well known and easy to prove (see for example \cite {GLM12}).\\
$(b)$ and $(c)$ are consequence of $(a)$.\hfill $\Box$.\smallskip

Boolean elements correspond with \emph{factor congruences}, in fact $\theta$ is a factor congruence of $\A\in\BRL$ if and only if there is $b\in B(\A)$ such that $\theta=\theta(\langle b\rangle)$, and so
 $\A$ is directly indecomposable if and only if $\B(\A)$ is the two element Boolean algebra (cf.\cite{KO01}). Since in congruence distributive varieties any  finitely subdirectly irreducible member is indecomposable, we have  that for any variety $\V\subseteq \BRL$, $\V_{\mbb{FSI}}\subseteq \V_{\mbb{DI}}$.\smallskip

 If  $F$ is an implicative filter of $\A$, then $F\cap B(\A)$ is a lattice filter of $\B(\A)$. An implicative filter  $F$ is called \emph{Stone filter} provided that $F=\langle F\cap B(\A)\rangle$, if,  in addition,   $F\cap B(\A)$ is ultrafilter, then $F$ is called \emph{Stone ultrafilter}. We represent by $St(\A)$ the set of all   Stone filters  of $\A$, and by $U(\A)$ the set of all Stone ultrafilters of $\A$.  

 \bg{lemma}\label{L:if-di} Let $F$ be an \imf\ of  $\A$. Then  $\A/F$ is directly indecomposable  if and only if it satisfies
   \bd
   \item[(di)]
   for any $a\in A$, $a\lor \neg a\in F$ implies $a\in F$ or $\neg a\in F$.
 \ed\e{lemma}
 {\em Proof:} $\Ra)$ Assume  $\A/F$  directly  indecomposable and $a\in A$  is such that $a\lor\neg a\in F$. Then $a/F\in B(\A/F)$, and so either $a/F=\top/F$ or $\neg a/F=\top/ F$, that is. $a\in F$ or $\neg a\in F$.

 \no $\La)$ If $F$ satisfies (di), and $a/F\in B(\A/F)$, then $a\lor \neg a \in F$, hence $a\in F$ or $\neg a\in F$, and so $a/F\in \{\top/F,\bot/F$\}. Therefore $\A/F$ is directly indecomposable.\hfill $\Box$

 \bg{theorem}\label{T:ufs-di} Let $F$ be an \imf\ of $\A$. Then
  \be[$(a)$]
  \item If $\A/F$ is directly indecomposable, then $\langle F\cap B(\A)\rangle\in U(\A)$
  \item If $F\in Sp(\A)$, then $A/F$ is directly indecomposable.
  \item If $F\in U(\A)$, then $A/F$ is directly indecomposable.
  \ee
 \e{theorem}
 {\em Proof:} $(a)$ and $(b)$ are consequence  of Lema  \ref{L:if-di}.
  To prove $(c)$ it is enough to see that  $F$ satisfies (di) whenever $F\in U(\A)$.
  Assume that $F\in U(\A)$ and $a\in A$ is such that  $a\lor \neg a\in F$. Then there is $b\in F\cap B(\A)$ such that $b\le a\lor \neg a$. Take $d=\neg b\lor a$, then, by $(a)$ and $(b)$ in Lemma \ref{L:previUltra}  we have
 \beqn
 d\lor \neg d &=& ( \neg b\lor a)\lor \neg(\neg b\lor a)=
 ( \neg b\lor a)\lor (b\land \neg a)\\
 &=&((b\land a)\lor \neg b)\lor (b\land \neg a)=((b\land a)\lor(b\land \neg a))\lor \neg b\\&=& (b \land(a\lor\neg a))\lor \neg b = b\lor\neg b=\top.\eeqn
 So $d\in B(\A)$, and since $F\cap B(\A)$ is an ultrafilter of $\B(\A)$, we have  $d\in F$ or $\neg d\in F$. Since $b\in F$, if   $d=\neg b\lor a= b\to a\in F$, then  $a\in F$ and if $\neg d= b\land\neg a\in F$, then $\neg a\in F$. Therefore $F$ satisfies (di).\hfill $\Box$\smallskip

\bg{theorem}\label{T:di=pri}  Let $\V$ be a variety of \brl s. Then the following are equivalent:
\be[$(i)$]
\item $\V_{\mbb{DI}}=\V_{\mbb{FSI}}$.
\item  for any $\A\in\V$, $U(\A)\subseteq Sp(\A)$.
\item for any $\A\in\V$, $U(\A)=Spm(\A)$.
\ee
\e{theorem}
{\em Prof:} $(i)$ implies $(ii)$: Assume that $\V_{\mbb{DI}}=\V_{\mbb{FSI}}$ and $\A\in \V$. If $F\in U(\A)$, then, by $(c)$ of Theorem \ref{T:ufs-di}, $\A/F\in \V_{\mbb{DI}}=\V_{\mbb{FSI}}$ and so  $F\in Sp(\A)$.

\no  $(ii)$ implies $(iii)$ because for any $P\in Sp(\A)$, $\langle P\cap B(\A)\rangle\in U(\A)$.

\no $(iii)$ implies $(i)$: If  $\A\in \V_{\mbb{DI}}$, then, by $(a)$ of Theorem \ref{T:ufs-di}, $\{\top\}\in U(\A)=Spm(\A)$ and so $\A\in \V_{\mbb{FSI}}$. This proves $\V_{\mbb{DI}}\subseteq\V_{\mbb{FSI}}$. The other inclusion always holds.\hfill$\Box$
\smallskip

Observe that, for each $P\in Sp(\A)$,  $\langle P\cap B(\A)\rangle$ is a Stone ultrafilter  contained in $P$,  and so   $\bigcap U(\A)\subseteq \bigcap Sp(\A)=\{\top\}.$
\smallskip

  Let $b(x)$ be an unary term and $\A\in \BRL$. We say that   \emph{$\A$ has $b(x)$ as Boolean term}, or that \emph{$b(x)$ is a Boolean term for $\A$},  whenever  for any $a\in A$ $b^{\A}(a)\in B(\A)$ and $b^{\A}(\top)\not= b^{\A}(\bot)$.  A variety  $\mbb{V}$  of \brls\ has  \emph{$b(x)$ as Boolean term} (or \emph{$b(x)$ is a Boolean term for $\V$}), whenever $b(x)$ is Boolean term for any $\A\in \V$.

  \bg{lemma}\label{L:bt in var} If $b(x)$ is a Boolean term for the variety $\V\subseteq \BRL$, then  there is $n>0$ such that  either $\V\models x^n\preccurlyeq b(x)$ or $\V\models x^n\preccurlyeq\neg b(x)$. In particular we have  either $\V\models b(\top)\approx \top$ or $\V\models \neg b(\top)\approx \top$.
  \e{lemma}
  {\em Proof:} Let $\F$ be the $1$-free algebra in $\V$ with $g$ as free generator.
  Since  $\langle g\rangle\in Max(\F)$ and $b^{\F}(g)\lor \neg b^{\F}(g)=\top\in \langle g\rangle$, then  either $b^{\F}(g)\in \langle g\rangle$ or $\neg b^{\F}(g)\in \langle g\rangle$. If $b^{\F}(g)\in \langle g\rangle$, then there is $n>0$ such that $g^n\le b^{\F}(g)$ and so $\V\models x^n\preccurlyeq b(x)$. In this case  $\top\le b^{\F}(\top)$, and so $\top= b^{\F}(\top)$.
  If $b^{\F}(g)\notin \langle g\rangle$, then $\neg b^{\F}(g)\in \langle g\rangle$, and the same argument shows that  there is $n>0$ such that $\V\models x^n\preccurlyeq\neg b(x)$. Thus  $\top\le \neg b^{\F}(\top)$, and so $\top= \neg b^{\F}(\top)$.
  \hfill $\Box$\smallskip

  It is easy to see  that for any $\A\in\BRL$, $B(\A)\cap Rad(\A)=\{\top\}$. Moreover,
  \bg{corollary}
  Let $b(x)$ be a Boolean term for $\A\in\BRL$ such that $b^{\A}(\top)\not=\bot$. Then:
  \be[$(a)$]
  \item There is $n>0$ such that for all $a\in A$ $a^n\le b^{\A}(a)$.
  \item  $Rad(\A)\subseteq\{a\in A: b^{\A}(a)=\top\}$.
  \ee
  \e{corollary}
  {\em Proof:} (a) Since $b(x)$ is a Boolean term  for the variety $HSP(\A)$ and $b^{\A}(\top)\not=\bot$, by Lemma \ref{L:bt in var}, there is $n>0$ such that $HSP(\A)\models x^n\to  b(x)=\top$.\\
  $(b)$ If $a\in Rad(\A)$ there is $k>0$, such that $\top=k.a^n\le k.b^{\A}(a)=b^{\A}(a)$. That shows $Rad(\A)\subseteq\{a\in A: b^{\A}(a)=\top\}$.\hfill $\Box$

  \begin{corollary}\label{C:kr->gap}
 Let $b(x)$ be an unary term. If all member of $\mbb{C}\subseteq\BRL$ has $b(x)$ as Boolean term, then $b(x)$ is Boolean term for the variety $HSP(\mbb{C})$.\hfill $\Box$
\end{corollary}

 Now we can characterize semisimple varieties by using Boolean terms,

  \bg{lemma}\label{SS-bterm} Let  $\V$ be variety of \brl s. Then $\V$ is  a semisimple variety if and only if
 there is $m>0$ such that  $x^m$ is  Boolean term for $\V$.
  \e{lemma}
  {\em Proof:} Assume that $\V$ is a semisimple variety, then, by Theorem \ref{T:PrevioSS}, there is $m>0$ such that $\V\subseteq \EMm$, and so $\V\models x\lor \neg x^m\approx \top$. Thus for any $\A\in \V$ and any $a\in A$, we have  $a\lor \neg a^m=\top$,  since $\langle a\rangle=\langle a^m\rangle$,
  \[
  \{\top\}= \langle a\lor \neg a^m\rangle=\langle a\rangle\cap \langle \neg a^m\rangle=
  \langle a^m\rangle\cap \langle \neg a^m\rangle =\langle a^m\lor \neg a^m\rangle.\]
  Therefore $a^m\lor \neg a^m= \top$ and hence $a^m\in B(\A)$, since $\top^n=\top\not=\bot=\bot^n$, the arbitrariness of $a$ and $\A$ proves that $x^m$ in a Boolean term for $\V$.

  \no To prove the converse  implication, observe that  for any $\A\in\V$ and any  $a\in A$  $a^m\le a$, and so $a^m\lor\neg a^m\le a\lor a^m$. Hence if $x^m$ is Boolean term for $\V$, then  $\V\models x\lor x^m\approx\top$ and so $\V\subseteq \EMm$.
  \hfill $\Box$

  \bg{remark}  It follows from Lemma \ref{SS-bterm} that for any $m>0$  $m.x^m$ is also a Boolean term for any variety contained in $\EMm$. However the converse is not true, because there are non-semisimple varieties of \brls\ having $m.x^m$ as Boolean term (cf. Theorem \ref{Gap<-> bool+rad}).
  \e{remark}
  
  In some varieties admit a Boolean term  that is also radical term. In this cases we have 
  
  \bg{corollary}\label{C:T=0->corad}  If a variety   $\V$ of \brl s admits $b(x)$ as boolean and radical term, then for any $\A\in\V_\mbb{DI}$, 
  \[A\smsm Rad(\A)=\{a\in A: b^{\A}(a)=\bot\}.\]  \e{corollary}

 Following \cite{CT12}, an unary term $t(x)$ is called \emph{Boolean retraction term for} $\A\in\BRL$ provided that the map $t^{\A}:a\mapsto t^{\A}(a)$ gives a retraction from $\A$ onto $\B(\A)$. \emph{A variety $\V\subseteq \BRL$   has Boolean retraction term},\label{P:BRT} if there is an unary term $t(x)$ such that it is Boolean retraction term for any member of $\V$.
 In fact, Boolean retraction term is a particular case of Boolean term. Varieties having  Boolean retraction term has been  studied in \cite{CT12}. Taking into account Theorem 2.7 of \cite{CT12}, we have

 \bg{lemma}\label{L:brt->bt rt} Let  $\V$ be a variety admitting $t(x)$ as Boolean retraction term. Then $t(x)$ is Boolean and radical term for $\V$.\hfill$\Box$
 \e{lemma}

Particular cases of unary terms are $m.x^n$, when $m,n>0$. In cases where one of these is a Boolean term for a certain variety, this one has some additional properties.

 \bg{theorem}\label{T:mnbt->mnradt} Let $n,m>0$ and let $\V$ a variety of \brl s admitting $m.x^n$ as Boolean term.   Then:
 \be[$(a)$]
 \item $\V\subseteq \WL_\mbb{k}$, for any $k\ge \max\{m,n\}$
 \item  $\V$ is $m$-radical
 \ee
 \e{theorem}
 \emph{Proof of ($a$)} Let $\A\in\V$. If   $a\in A$ and $k\ge\max\{m,n\}$, then  $a^{k}\le a^n\le m.a^n\le k.a^n\le k.a$, and so
 $\top =m.a^n\lor \neg m.a^n
 \le k.a\lor\neg a^k =k.a\lor k.\neg a.$. Hence $\A\in\WL_\mbb{k}$..\\
   \emph{Proof of $(b)$} By \cite[Lemma 1.8]{CT12} and $(a)$, $\V$ is $k$-radical. Moreover, if  $\A\in V$ and $a\in Rad(\A)$, then for any $s\le m$, $\top=m.a^n\le m.a^s$, and for any  $s>n$, there is $r>0$  such that $rn<s\le (r+1)n$ and so, since $a^r$ and $a^{r+1}$ are in $ Rad(\A)$,
   $\top = m.(a^{r+1})^n\le m.a^s\le m.(a^r)^n=\top$. Thus $Rad(\A)\subseteq \{a\in A: \forall s>0\ m.a^s=\top\}$. The other inclusion is always  true. \hfill $\Box$

\section{Local algebras and Generalized Apple\\ Pro\-perty}

An algebra $\A$ is called \emph{local} provided that it has a unique proper maximal congruence, i.e, $\langle Con(\A),\subseteq \rangle$ has  only one coatom. 
 \bg{lemma}\label{L:DI_is_L} If $\A$ is local algebra, then $\A$ is directly indecomposable.
 \e{lemma}
 \emph{ Proof:}
 Assume  $\A$ is a local. If $\A$ is not    directly indecomposable, then there is a non-trivial pair factor of congruences  $\{\theta_1,\theta_2\}$. For $i\in\{1,2\}$, $\theta_i\not=\nabla_A$, then, since $\langle Con(\A),\subseteq \rangle$ is closed under the union of upward directed families,  by Zorn's Lemma,   there is a proper maximal congruence $\varphi_i$ containing $\theta_i$. Moreover,   since $ \theta_1\vee\theta_2=A$, we have  $\varphi_1\not=\varphi_2$. This contradicts the assumption, and  so   $\A$ is directly indecomposable.\hfill$\Box$ \smallskip

The following lemma is easily provable.
\bg{lemma}\label{L:GLocal}
 Let  $\theta$ be a proper   congruence relation  on an algebra $\A$, then it satisfies:
 \be[$(a)$]
 \item   $\A/\theta$ is local if and only if $\theta$ is contained in  only one proper maximal congruence  of $\A$.
 \item If $\A$ is local, then $\A/\theta$ is local.\hfill $\Box$
 \ee
 \e{lemma}

 In what follows, given a class of algebras $\mbb{K}$, $\mbb{K_L}$ denotes the local members of $\mbb{K}$. Then, by the above lemma,
 for any variety $\V$, $\V_{\mbb{L}}\subseteq \V_{\mbb{DI}}$. A variety $\V$ is called \emph{locally representable}  if all its members  are representable as subdirect product of local algebras. This is equivalent to $\V_{\mbb{SI}}\subseteq \V_{\mbb{L}}$.

 \bg{remark}\label{R:WLK-LR}It is shown in \cite[Lemma 4.6]{T16} that for any $k>0$ $\WLk_{\mbb{FSI}}\subseteq \WLk_{\mbb{L}}$, and so  any subvariety of $\WLk$ is locally representable. Therefore any prime i-filter of any member of $\WLk$ is contained in only one maximal i-filter (\cite[Corollary 4.7]{T16})
 \e{remark}

 We say that an algebra $\A$ has the {\em General Appel Property} (GAP for short),   if   $\theta\in Con(\A)\smsm\{\nabla_A\}$ is such that $\A/\theta$ is directly indecomposable, then   the interval lattice $[\,\theta,\nabla_A]$ has a coatom, i.e., $\A/\theta$ is local. And we say that a variety $\V$ has the \emph{General Apple Property} (GAP) provided that any $\A\in \V$ has GAP. Then from Lemmas \ref{L:DI_is_L} and \ref{L:GLocal} we deduce the following
 \bg{lemma}\label{L:GAP-Locrep} A variety $\V$ has GAP if and only if $\mbb{V_{DI}}=\mbb{V_L}$. \hfill$\Box$
 \e{lemma}

 \begin{remark}\label{R:finite GAP} Observe that a finite algebra having GAP also has the Apple property, and so, any variety with GAP also has the Apple property.
 \end{remark}

  We return to \brl . In what follows $\A$ represents a \brl\ and $\V$ a subvariety of $\BRL$. From the properties described in Section 2, we deduce that  $\A$ has the GAP provided that for any implicative filter $F$ of $\A$ if $\A/F$ is directly indecomposable, then $F$    is contained in a unique maximal implicative filter. We can improve this equivalence, because   from Theorem \ref{T:ufs-di} it follows:
 \begin{lemma}\label{GAP-Stone}
 A \brl\ $\A$ has  GAP if and only if every Stone ultrafilter of $\A$ is contained in a unique maximal implicative filter.\hfill$\Box$
 \end{lemma}

 As a consequence we have
 \bg{lemma} Let $\A$ be a \brl\ having  GAP. Then it sa\-tisfies:
 \be[$(a)$]
 \item Any prime implicative filter of $\A$ is contained in a unique maximal implicative filter.
 \item If $P,P'\in Sp(\A)$ are such that $P\cap B(\A)=P'\cap B(\A)$, then $P$ and $P'$ are contained in the same maximal implicative filter.\hfill$\Box$
 \ee
 \e{lemma}

Then   (cf. \cite[Theorem 3.3]{CT96}, \cite[Proposition 4.5]{DGL00} and \cite[Theorem 9]{L05})
 \bg{theorem}\label{T:rad(f)-max} For any $\A\in\BRL$ the following properties  are equivalent:
 \be[$(i)$]
 \item $\A$ has  GAP
 \item For any $a\in A$ there is $b\in B(\A)$ such that $b\in\langle a\rangle$ and $\neg b\in \langle\neg a\rangle$
 \item For any $a \in A$ there are $b \in B(A)$ and $m, n > 0$, such that\\ $a^m \le b \le n.a$.
 \item For any $a\in A$ there are $b\in B(\A)$ and $k>0$, such that $a^k\le b\le k.a$ \ee
\e{theorem}
{\em Proof.}
 $(i)\Ra (ii)$: Suppose $a\in A$, and let $G=\langle a\rangle\cap B(\A)$ and $F=\{b\in B(\A): \neg b\in \langle \neg a\rangle\}$. We need to see that  $G\cap F\not=\emptyset.$ Suppose not, i.e., $G\cap F=\emptyset$.
 It is clear that $G$ and $F$ are respectively a filter and an ideal of $\B(\A)$, and hence there is $P\in Sp(\B(\A))$ such that $G\subseteq P$ and $P\cap F=\emptyset$. Then we have that $\neg a\notin \langle \{a\}\cup P\rangle$ and $a\notin \langle \{\neg a\}\cup P\rangle$.
 Indeed, if $\neg a\in \langle \{a\}\cup P\rangle$, then there would be $n>0$  and $b\in P$ such that $b\cdot a^n\le \neg a$.
 Hence $b\cdot a^{n+1}\le a\cdot \neg a =\bot$, and so $a^{n+1}\le \neg b$, i.e., $\neg b\in G\subseteq P$, and $P$ would not be  proper.
 If $a\in \langle \{\neg a\}\cup P\rangle$, then there would be $n>0$ and $b\in P$  such that $b\cdot (\neg a)^n\le  a$. Then $ b\cdot(\neg a)^{n+1}\le a\cdot\neg a=\bot$,  so $(\neg a)^{n+1}\le \neg b$, and would have $b\in P\cap F=\emptyset$. Therefore $\langle \{a\}\cup P\rangle$ and $\langle \{\neg a\}\cup P\rangle$ are proper \imf s, and hence there are maximal implicative filters  $M_1$ and $M_2$ such that $\langle \{a\}\cup P\rangle\subseteq M_1$ and $\langle \{\neg a\}\cup P\rangle\subseteq M_2$. Since $\langle \{a,\neg a\}\rangle =A$. $\neg a\notin M_1$ and $a\notin M_2$, and $M_1\not=M_2$. But  the Stone ultrafilter $\langle P\rangle$ is contained in $ M_1\cap M_2$, that contradicts Lemma \ref{GAP-Stone}. Therefore $F\cap G\not=\emptyset$.

 \no $(ii)\Ra (i)$.  Suppose $(ii)$ holds and let $M_1$ and $M_2\in Max(\A)$  distinct maximal \imf s of $\A$. Let $a\in M_1$, and $a\not\in M_2$. There is $n>0$ such that $\neg a^n\in M_2$. Then there is $b\in B(\A)$ such that $b\in \langle a^n\rangle=\langle a\rangle\subseteq M_1$ and $\neg b\in \langle \neg a^n\rangle\subseteq M_2$. Then $\langle M_1\cap B(\A)\rangle\not = \langle M_2\cap B(\A)\rangle$. Hence there is no Stone ultraideals contained in $M_1\cap M_2$.

 \no Assume  $(ii)$  and consider  $a\in A$, then there exist $b\in B(\A)$ and $m,n>0$ such that $a^m\le b$ and $(\neg a)^n\le \neg b$, and so $a^m\le b$ and $b=\neg\neg b\le \neg(\neg a)^n=n.a$. This proves $(ii)$ implies $(iii)$

 \no To see that $(iii)$ implies $(iv)$ it suffices take
  $k=\max\{n,m\}$, and apply the properties stated in Lemma \ref{L:nmplus}.

 \no Finally, $(iv)$ implies $(ii)$ follows from the fact that for any $a\in A$ and any $b\in B(\A)$ and any $k>0$, $(\neg x)^k\le \neg b$ if and only if $b\le k.a$.\hfill $\Box$\smallskip

 The  Property $(iii)$,   so  $(i)$,  $(ii)$ and $(iv)$ of Theorem \ref{T:rad(f)-max}  characterizes the \brls\ having the \emph{lifting Boolean property} (BLP), i.e., for any implicative filter $F$, $B(\A/F)= B(\A)/F$(see \cite[Proposition 4.5]{GM14}). These \brls\ are also called \emph{quasi-local}. We add a simple proof of this fact.

 \begin{lemma}[cf,\mbox{\cite[Propositions 4.5 and 4.16]{GM14}}]
 Let $\A$ be a \brl\ $\A$. Then $\A$ has BLP if and only if $\A$ satisfies $(iv)$ of Theorem \ref{T:rad(f)-max} if and only if it has GAP.
 \end{lemma}
 {\em Proof:} $ \Ra)$ Assume  $\A$ has BLP  and  $a$ is an arbitrary element of $ A$. If $F= \langle a\lor \neg a\rangle$, then $a/F\in B(\A/F)=B(\A)/F$. Hence there is $b\in B(\A)$ such that $b/F=a/F$ and $\neg a/F=\neg b/F$, and  so $a\to b\in F$ and $\neg a\to \neg b\in F$. Therefore exists $k>0$ such that $(a\lor\neg a)^k\to (a\to b)=\top$ and $(a\lor\neg a)^k\to (\neg a\to \neg b)=\top$.
 Thus
 \begin{align*}
 \top &=(a\lor\neg a)^k\to (a\to b)= (a^k\lor (\neg a)^k)\to (a\to b)\\
 &=(a^k\to (a\to b))\land ( (\neg a)^k\to (a\to b))
 = a^{k+1}\to  b,
 \end{align*}
 that is,  $a^{k+1}\le  b$. Similarly, from $(a\lor\neg a)^k\to (\neg a\to \neg b)=\top$  we would prove  $(\neg a)^{k+1}\le \neg b$, and so $b=\neg\neg b\le \neg(\neg a)^k\le(k+1).a$.

 \no$\La)$ Assume  $\A$ satisfies $(iii)$ of \ref{T:rad(f)-max}. Take $F$  and \imf\
 and $a\in A$. Then there are $b\in B(\A)$ and $k \ge 0$  such that $a^k\le b\le k.a$. Take $F$ an \imf\   such that $a/F\in B(A/F)$, then
 \[a/F=(a/F)^k=a^k/F\le  b/F\le k.a/F =k.(a/F)=a/F.\]
  Thus  $a/F=b/F\in B(\A)/F$. This shows  that $B(A/F)\subseteq B(\A)/F$. The reverse inclusion is always true.\hfill $\Box$\smallskip

 The following theorem  gives an equational characterization of varieties of \brl s having GAP.
 \bg{theorem}\label{T:GAP_bt} For any variety of $\V \subseteq \BRL$ the following conditions are equivalent:
 \be[$(i)$]
 \item  $\V$ has   GAP.
  \item There exist $b(x)$  Boolean term for $\V$ and $n,m>0$ such that
     \bg{equation}\V\models (x^m\to b(x))\cdot(b(x)\to n.x)\approx \top.\label{Eq:VDI=VL1}
     \e{equation}
 \item There exist $b(x)$  Boolean term for $\V$ and $k>0$ such that
     \bg{equation}\V\models (x^k\to b(x))\cdot(b(x)\to k.x)\approx \top.\label{Eq:VDI=VL}
     \e{equation}
 \ee
 \e{theorem}
 {\em Proof:}  $(i)$ implies $(ii)$.  Let $\F$ be the 1-free algebra of $\V$ with free generator $g$. By Theorem \ref{T:rad(f)-max}, there are an unary term $b(x)$ and $k>0$ such that   $b^{\F}(g)\in B(\F)$ and $g^k\le b^{\F}(g)\le k.g$. Therefore, $b^{\F}(g)\lor \neg b^{\F}(g)=\top$, $g^k\to b^{\F}(g)=\top$ and $ b^{\F}(g)\to \neg(\neg g^k)=\top$.  Hence $b(x)$ is a Boolean term for $\V$ (see Lemma \ref{L:bt in var}), and taking into account item $(\ref{bb8})$ in Lemma \ref{basicbounded},   $(\ref{Eq:VDI=VL})$ holds.

 \no To see  $(iii)$ implies $(i)$ it suffices take $k=\max\{m,n\}$

 \no $(iii)$ implies $(i)$ follows from Theorem~\ref{T:rad(f)-max} and Lemma \ref{L:bt in var}, using the pro\-perties listed in Lemma \ref{basicbounded}.
 \hfill $\Box$

\bg{remark} It follows from the above result that a variety $\V\subseteq\BRL$ has GAP if and only if there exists an unary term $b(x)$ and there exists $k>0$ such that  $\V$ satisfies the equation $b(x)\lor\neg b(x)\approx\top$ and $(\ref{Eq:VDI=VL})$.  Since $\V\models b(\top)= \top$ and $\V\models b(\bot)=\bot$,  $b(x)$ is Boolean term for $\V$.
\e{remark}
Since for any $k>0$ $\BRL\models x^k\preccurlyeq k.x^k\preccurlyeq k.x$, we have
\bg{corollary}\label{C:boole-nm} If there is  $k>0$ such that  $k.x^k$ is Boolean term for  a variety $\V\subseteq \BRL$, then $\V$ and $b(x)=k.x^k$ satisfy (\ref{Eq:VDI=VL}) of Theorem \ref{T:GAP_bt}, and so $\V$ has GAP.\hfill $\Box$
\e{corollary}
 \bg{lemma} Assume $b(x)$ is a Boolean term  for the variety $\V\subseteq \BRL$ such that  (\ref{Eq:VDI=VL}) is satisfied for  $k>0$. Then the following properties hold true:
 \be[$(a)$]
 \item $\V\models b(x^k)\preccurlyeq b(x)$.
 \item $\V\models m.b(x)\approx b(x)$, for any $m>0$.
 \item $\V\models b(b(x))\approx b(x)$.
 \item $\V\models x^k\preccurlyeq \neg b(\neg x)\preccurlyeq k.x$.
 \item  $\V\models k.x^k\preccurlyeq (k.x)^k$
\ee
 \e{lemma}
 \emph{Proof:} To prove these properties we consider $a$ an arbitrary element of $\A\in\V$.\\
 $(a)$: Since $a^k\le b^{\A}(a)$ and $b^{\A}(a)\in B(\A)$, we have \[ b^{\A}(a^k)\le k.a^k\le k.b^{\A}(a)= b^{\A}(a).\]
  $(b)$ follows from Lemma \ref{L:previUltra}.\\
 $(c)$: Since $b^{\A}(a)\in B(\A)$,we have
 \[b^{\A}(a)= b^{\A}(a)^k\le b( b^{\A}(a))\le k.b^{\A}(a)=b^{\A}(a).\]
 $(d)$: By (\ref{Eq:VDI=VL})
 $(\neg a)^k\le  b^{\A}(\neg a)\le k.\neg a=\neg a^k$, and so
 \[a^k\le \neg\neg a^k\le \neg b^{\A}(\neg a)\le \neg (\neg a)^k=k.a.\]
$(e)$: since $b^{\A}(a)\in B(\A)$,  we have
 \[k.a^k\le k.b^{\A}(a) = b^{\A}(a)= b^{\A}(a)^k\le  (k.a)^k\]

 \begin{corollary}\label{C:kr->gap} Let $b(x)$ be an unary term.
 If all member of $\mbb{C}\subseteq\BRL$ has $b(x)$ as Boolean term satisfying (\ref{Eq:VDI=VL1}) then $b(x)$ is Boolean term for the variety $HSP(\mbb{C})$ and so it has GAP.\hfill $\Box$
 \end{corollary}

 As a consequence of Theorem \ref{T:GAP_bt} we have:
 \bg{lemma}\label{L:Gap->wlk}
 If the variety $\V\subseteq \BRL$ has  GAP, then  there is $k>0$ such that $\V\subseteq \WLk$.
 \e{lemma}
 {\em Proof:} Assume $\V$ has the GAP, then there is a Boolean term $b(x)$ and $k>0$ such that $(\ref{Eq:VDI=VL})$ holds in $\V$. Let $\A\in\V$, for   any $a\in A$ we have  $a^k\le b^{\A}(a)\le k.a$, and so
 \[\top =b^{\A}(a)\lor \neg b^{\A}(a)= \neg\neg b^{\A}(a)\lor \neg b^{\A}(a)
 \le \neg k.a \lor a^k =k.a\lor k.\neg a.\]
 And $\A\in \WLk$.\hfill $\Box$

 \begin{remark} Observe that if  $\A\in \BRL$  has a Boolean term satisfying $(3.2)$ for $k>0$, then it is possible that $\A\in \WL_{\mbb r}$, for some $0<r<k$, as seen in  Example \ref{Ex:1} detailed in the Appendix.\end{remark}
 
 Let $r,k$ be non-negative integers. From the properties described in  Lemma \ref{L:nmplus} it follows that for each $m\ge r,k$ $\BRL\models x^m\preccurlyeq k.x^r$ and $\BRL\models k.x^r \preccurlyeq m.x$. In view of Corollary \ref{C:wlk-trad}, we have

 \bg{theorem}[cf. Theorem \ref{T:Qe-k-rad}]\label{Gap<-> bool+rad}Let $\V$ be a subvariety of $\BRL$. Then the following conditions  are equivalent
 \be[$(i)$]
 \item There is $k>0$  such that $\V\subseteq \WL_k$ and $\V_{\mbb{SS}}$ is a variety.
 \item There is $k>0$  such that $\V\subseteq \WL_k$ and $\V$ admits radical term.
 \item  There is $k,r>0$ such that $k.x^{kr}$  is  Boolean and radical term for $\V$.
 \item  $\V$ admits Boolean   and radical term satisfying (\ref{Eq:VDI=VL1}).
 \item $\V$ has GAP and admits radical term.

   \ee
  \e{theorem}
   {\em Proof:} $(i)$ iff $(ii)$ follows from Corollary \ref{C:wlk-trad}.

  \no $(ii)$ implies $(iii)$. Assume that there is $k>0$ such that  $\V\subseteq \WLk$, and $\V$  admits radical term.  Then,   by Corollary \ref{C:wlk-trad}, there exists   $r\ge k$ such that  $k.x^r$ is radical term for $\V$, so (see  Remark \ref{R:krads>r}) $k.x^{kr}$ is also radical term for $\V$. We are going to see that $k.x^{kr}$ is  Boolean term for $\V$.
   Take $\A\in \V_\mbb{FSI}$, then, by Remark \ref{R:WLK-LR}, $\A\in\V_{\mbb{L}}$. Let  $a$ be an arbitrary element of $A$, such that $k.a^{kr}\not=\top$. Then $a\notin Rad(\A)$ and $k.a^r\not=\top$, and since $\A\in\WLk_\mbb{FSI}\subseteq\WLk_\mbb{DI}$, we have $\top =k.\neg a^r =\neg a^{rk}$. Therefore  $a^{rk}=\bot$, and $k.a^{rk}=\bot$. That shows that $k.x^{rk}$ is  Boolean term and radical term for any $\A$ in $\V_{\mbb{FSI}}$. Now since each member of $\V$ is subdirect product of members of $\mbb{V_{FSI}}$, by Lemmas \ref{L: SbProdH} and \ref{L: SbProdH2}, $k.x^{rk}$ is also Boolean  term  for $\V$.

    \no $(iii)$ implies $(iv)$, is trivial.

      \no $(iv)$ implies $(v)$, because by  Corollary \ref{C:boole-nm}, $(iv)$ implies  $\V$ has the GAP.

   \no  $(v)$ implies $(ii)$, because   if $\V$ has GAP, then $\V\subseteq\WLk$ for some $k>0$.
   \hfill $\Box$\smallskip
   
   The following lemma shows that Boolean  and radical term  satisfying (\ref{Eq:VDI=VL1}), or (\ref{Eq:VDI=VL}), are essentially unique on varieties of bounded residuated lattices.
   
 \bg{lemma} If a variety $\V\subseteq \BRL$ admits $s(x)$ and $t(x)$ as Boolean and radical terms satisfying (\ref{Eq:VDI=VL1}), then the equation $s(x) \approx t(x)$ holds in $\V$.
   \e{lemma}
 \emph{Proof:}  Since $\V$ has GAP, It suffices that the equation holds in any $\A\in \V_\mbb{DI}$. By Corollary, \ref{C:T=0->corad} if  $a\in A$, then 
 \begin{eqnarray*}
    s^{\A}(a)=\top & \mbox{iff} & a\in Rad(\A)\mbox{ iff }t^{\A}(a)=\top \\
    s^{\A}(a)=\bot & \mbox{iff} & a\notin Rad(\A)\mbox{ iff }t^{\A}(a)=\bot 
  \end{eqnarray*}
  An so $s^{\A}(a)=t^{\A}(a)$.\hfill $\Box$\smallskip

 By Lemma \ref{L:brt->bt rt}, \emph{ if a variety $\V$ has $t(x)$ as Boolean retraction term, then   $t(x)$ is also radical term for $\V$ and $\V_{\mbb{SS}}$ is the variety of Boolean algebras}; thus   there are $k,r>0$ such that  $k.x^r$  is a radical term for $\V$, and so for any $n\ge r$  $k.x^n$ is also radical term for $\V$. In \cite[Section 5]{CT12} it is shown that  for any $k>0$, the subvarieties of $\BRL$ defined by the equation $k.x^k\approx (k.x)^k$, denoted by $\Vk$, is the greatest subvariety of $\WLk$ having a Boolean retraction term. These facts allow us to improve theorems 5.6 and 5.7 of \cite{CT12}, showing that a variety $\V\subseteq \BRL$ has Boolean retraction term if and only if there is $k>0$ such that $\V\subseteq \Vk$.

 \bg{theorem}[cf.Theorem 5.7 of \cite{CT12}]\label{T:Br-WLn}For any subvariety $\V$ of $\BRL$, the following conditions are equi\-valent:
 \be[$(i)$]
 \item $\V$ admits Boolean retraction term.
 \item $\V\subseteq \WLk$ and $\V\models k.x^k\approx(k.x)^k$, for some $k>0$.
 \item $\V$ satisfies the equation $k.x^k\lor k.(\neg x)^k\approx\top$, for some $k>0$.
 \item There exists $k>0$ such that $k.x^k$ is Boolean  retraction term for $\V$.
 \ee
 \e{theorem}
 {\em Proof:}    $(i)$ implies $(ii)$: If $\V$ has Boolean retraction term, then, it admits Boolean and radical term. By Theorem \ref{Gap<-> bool+rad}, there is $k>0$ such that $\V\subseteq \WLk$ and, by \cite[Theorem 5.6]{CT12}, $\V\models k.x^k\approx(k.x)^k$.

  \no $(ii)$ implies $(iii)$ follows from \cite[Theorem 5.7]{CT12}.
    
  \no $(iii)$ implies $(iv)$ follows from \cite[Theorems 5.1 and 5.6]{CT12}. 
    
  \no $(iv)$ implies $(i)$ is evident.\hfill $\Box$
  \smallskip

Recall that for each $m>0$ $m.x^m$ is radical  term for $\Em$.
\bg{lemma} \label{L:sEm-loc-repr}Let $\A\in\Em_{\mbb{DI}}$.  Then $\A$ is local if and only if $m.x^m$ is Boolean term for $\A$.
\e{lemma}
{\em Proof:}  Assume that $\A$ is local. Let $a\in A$, since  $Rad(\A)$ is maximal i-filter and $\A/Rad(\A)\in \EMm$, $a\lor \neg a^m\in Rad(\A)$. If $a\in Rad(\A)$, then $m.a^m=\top$. If $a\notin Rad(\A)$, then  $\neg a^m\in Rad(\A)$,   and  by  item $\ref{plus61})$ of  Lemma \ref{L:nmplus} we have
\[\top = (m.(\neg(a^m))^m\le m.\neg a^m= \neg a^{mm}\le \neg a^m,\]
hence  $a^m =\bot$ and so $ m.a^m=\bot$. Therfore $m.x^m$ is Boolean term for $\A$.\smallskip

 \no If $m.x^m$ is Boolean term for $\V$, then, since  for any $a\in A$, $a^m\le m.a^m\le m.a$, by Corollaries  \ref{C:boole-nm} and \ref{C:kr->gap},   $HSP(\A)$ has  GAP.  Then, since $\A$ is directly indecomposable,  $\A$ is local. \hfill$\Box$
 \smallskip

Any variety is generated by is directly indecomposable members, so
\bg{theorem} \label{T:loca,boterm} Let $\V$ a subvariety of $\Em$. Then the following are equivalent
 \be[$(i)$]
  \item $\V$ has GAP
  \item $m.x^m$ is a Boolean term for $\V$.
  \item $\V\subseteq\WLk$, for some $k>0$. 
  \ee
\e{theorem}

We have several immediate consequences from the above result
\begin{corollary}\label{C:loca,boterm1}
  A finite algebra in $\WLk$ has the GAP, and so AP.
\end{corollary}
{\em Proof:} It deduces from the above theorem, because any finite algebra in $\WLk$ belongs to $\Em$, for some $m>0$. By Remark \ref{R:finite GAP} it has AP. \hfill$\Box$ \medskip

Taking into account Corollary \ref{C:lf->m-cont}, we  deduce
\bg{corollary} Any locally finite variety of \brls\ has GAP.\hfill $\Box$
 \e{corollary}

 \appendix
 \section{Two examples}
 To end the paper we  present two examples. However, an infinite denume\-rable family of involutive \brl\ having Boolean terms satisfying $(\ref{Eq:VDI=VL1})$ can be found in \cite{T23}, which is currently being drafted.
 
 \bg{example}\label{Ex:1}{\rm We assume that the reader  is familiar  with  MV-algebras and Wajsberg hoops. (see  for example \cite{CDM00}, \cite{GMT99} and \cite{BF00}  and the references given there).

Let $\langle\mathds{Z};+,0;\le\rangle$ be the totally ordered group of integers, and let $\mathds{Z}\stackrel{\to}{\otimes}\mathds{Z}=\langle \mathds{Z}\times\mathds{Z},+,(0,0);\preccurlyeq  \rangle$ be  the lexicographic product  whose order $\preccurlyeq$ is defined:
\[(a,r)\preccurlyeq (b,s) \mbox{ iff }\left\{\begin{array}{ll}
 a<b, & \mbox{ or   } \\
  a=b & \mbox{ and }r\le s\end{array}\right.\]

 \no For any  integer $n>0$ we consider the MV-chain $\bfm{\LLn}=\bfm{\Gamma}(\mathds{Z}\stackrel{\to}{\otimes}\mathds{Z},(n,0))$, obtained from the totally ordered grup $ \mathds{Z}\stackrel{\to}{\otimes}\mathds{Z}$ with the strong unit $(n,0)$, see \cite{CDM00}, considered as  a  involutive \brl\ (Wajsberg algebra or bounded Wajsberg Hoop), i.e.,
 \bi\item[]$\bfm{\LLn}=\left\langle L^{\omega}_{n+1}=[(0,0),(n,0)];\ast,\to,\land, \lor,\bot=(0,0),\top=(n,0)\right\rangle$;\ei
 such that
 \bi
 \item its universe is the interval
  \[[(0,0),(n,0)]=\{(a,r)\in \mathds{Z}\times\mathds{Z}:(0,0)\preccurlyeq(a,r)\preccurlyeq(n,0) \},\]
  \item its lattice order $\le_{\bfm{\LLn}}= \ \preccurlyeq\! |_{\LLn}$ is total,
  \item $(a,r)\to (a,s)=\min\{(n,0),[(n,0)-(a,r)]+(b,s)\}$\footnote{``$\max$" and ``$\min$" are relative  to order $\preccurlyeq$ of $\mathds{Z}\stackrel{\to}{\otimes}\mathds{Z}$},
  \item $(a,r)\ast (b,s)= \max\{(0,0),[(a,r)+ (b,s)]-(n,0)\}=\neg ((a,r)\to\neg (b,s)),$\\
  where $\neg x=x\to \bot$.
      \ei

  We also consider $L_{n+1}=\{\bfm{m}=(m,0)\in \mathds{Z}\times \mathds{Z}: 0\le m\le n\}$ which is universe of $\bfm{L_{n+1}}$  subalgebra of $\bfm{\LLn}$.  $\bfm{L_{n+1}}$   is also a copy of the  totally ordered Wajsberg hoop
  with $n+1$ elements with lower bound (see \cite{BF00}). Moreover, $\bfm{L_{n+1}}$ is simple. In fact, any finite simple MV-algebra is isomorphic to   $\bfm{L_{n+1}}$, for some $n>0$.
  In par\-ticular, $\bfm{L_{2}}$ is a copy of the two element Boolean algebra, in which  $\ast=\land$ and $\neg x $ is the complement of $x$. 
  It is well known and easy to chek that the following properties are satisfied for all $n>0$,
  
  \bi
  \item If $m>n$, then $2.x^{m}$ is Boolean and radical term for $HSP(\bfm{L_{n+1}^\omega})$.
  \item $2.x^{n}$ is not Boolean term for $HSP(\bfm{L_{n+1}^\omega})$, because
  \item[]  \qquad  $2.(n-1,1)^n=(0,2n)\not=\bot$, $2.x^{n}$.
  \item $HSP(\bfm{L_{n+1}^\omega})_{\mbb{SS}}=HSP(\bfm{L_{n+1}})\subseteq \EMn$
  \ei
  Observe that by Lemma \ref{L:Gap->wlk} $\bfm{L_{n+1}^\omega}\subseteq \mbb{WL_\mbb{n+1}}$, although as is well known  $\bfm{L_{n+1}^\omega}\in \mbb{MV}\subseteq \mbb{MTL}\subseteq \mbb{WL_\mbb{2}}$.}
  \e{example}\newpage
 
 \bg{example}\label{Ex:2}{\em Consider the algebra
 $\bfm{W_9}=\langle W_9;\cdot,\to, \land,\lor,\bot,\top\rangle\}$ of type $(2,2,2,2,0,0)$
 such that $W_9=\{\top,1,2,3,4,5,6,7,\bot\}$, $\langle W_9; \land,\lor,\bot,\top\rangle$ is the bounded  lattice  given by the diagram  depicted in  Figure W9 and the operations $\cdot$ and $\to $ are given in the next tables:

 ${}$\vspace{1cm}

 \bg{figure}[ht]\bg{center}{\scriptsize
  \bg{picture}(0,100)
  \put(-180,100){\begin{tabular}{c||c|c|c|c|c|c|c|c|c}
   $\cdot$ & $\top$&$1$ &$2$ &$3$ &$4$ &$5$ &$6$ &$7$ &$\bot$ \\\hline \hline
   $\top$ &$\top$&$1$ &$2$ &$3$ &$4$ &$5$ &$6$ &$7$ &$\bot$ \\  \hline
      1   &  1   &$3$ &$3$ &$3$ &$5$ &  5  & 7  &$\bot$&$\bot$ \\\hline
   $  2 $ &  2   & 3  & 3  & 3  & 5  &  5  & 7  &$\bot$&$\bot$  \\ \hline
      3   &  3   & 3  & 3  & 3  & 5  &  5  &$\bot$&$\bot$&$\bot$  \\  \hline
   $  4 $ &  4   & 5  & 5  & 5  & 7  &$\bot$& 7   &$\bot$&$\bot$ \\ \hline
   $  5 $ &  5   & 5  & 5  & 5  &$\bot$&$\bot$&$\bot$&$\bot$&$\bot$  \\ \hline
      6   &  6   & 7  & 7  &$\bot$& 7  &$\bot$&$\bot$&$\bot$&$\bot$  \\ \hline
   $  7 $ &  7&$\bot$&$\bot$&$\bot$&$\bot$&$\bot$&$\bot$&$\bot$&$\bot$  \\ \hline
   $\bot$ &$\bot$&$\bot$&$\bot$&$\bot$&$\bot$&$\bot$&$\bot$&$\bot$ &$\bot$  \\
 \end{tabular}}
  \put(-180,-10){\begin{tabular}{c||c|c|c|c|c|c|c|c|c}
   $\to$ & $\top$&$1$ &$2$ &$3$ &$4$ &$5$ &$6$ &$7$ &$\bot$ \\\hline \hline
   $\top$ &$\top$&$1$ &$2$ &$3$ &$4$ &$5$ &$6$ &$7$ &$\bot$ \\  \hline
      1   &$\top$&$\top$& 1     & 1  & 4  & 4  & 6 & 6  & 7   \\\hline
   $  2 $ &$\top$&$\top$&$\top$ & 1  & 4  & 4  & 6 & 6  & 7   \\ \hline
      3   &$\top$&$\top$&$\top$&$\top$& 4 & 4  & 6 & 6 & 6  \\  \hline
   $  4 $ &$\top$&$\top$& 1  &1&$\top$&$1$& 4  & 4 & 5  \\ \hline
   $  5 $ &$\top$&$\top$&$\top$&$\top$&$\top$&$\top$&4& 4 & 4  \\ \hline
      6   &$\top$&$\top$& 1  & 1&$\top$& 1 &$\top$&$ 1 $&3 \\ \hline
   $  7 $ &$\top$&$\top$&$\top$&$\top$&$\top$&$\top$&$\top$&$\top$&1  \\ \hline
   $\bot$ &$\top$&$\top$&$\top$&$\top$&$\top$&$\top$&$\top$&$\top$& $\top$ \\
 \end{tabular}}
 \put(70,40){\circle*{3}}
 \put(100,70){\circle*{3}}
 \put(55,55){\circle*{3}}
 \put(70,70){\circle*{3}}
 \put(85,25){\circle*{3}}
 \put(115,55){\circle*{3}}
  \put(85,85){\circle*{3}}
 \multiput(85,25)(-15,15){3}{\line(1,1){30}}
 \multiput(85,25)(30,30){2}{\line(-1,1){30}}
 \put(85,85){\line(0,1){25}}
 \put(85,0){\line(0,1){25}}
 \put(85,0){\circle*{3}}
 \put(85,110){\circle*{3}}
 \put(85,-10){\makebox(0,0){$\bot$}}
 \put(85,120){\makebox(0,0){$\top$}}
  \put(80,85){\makebox(0,10){$1$}}
  \put(65,75){\makebox(0,0){$2$}}
 \put(48,55){\makebox(0,0){$3$}}
  \put(105,75){\makebox(0,0){$4$}}
 \put(62,35){\makebox(0,0){$5$}}
\put(120,55){\makebox(0,0){$6$}}
 \put(75,20){\makebox(0,0){$7$}}
 \put(85,-20){\makebox(0,0){Figure W9 : Lattice reduct of $\bfm{W_9}$}}
  \e{picture} }
 \end{center}
\end{figure}
${}$\vspace{1cm}

It is easy to verify that $\bfm{W_9}$ is a   \brl\ such that:
\bi
\item Its monoidal operation is $\cdot$, and if $\le$ is its lattice partial order, then   for any $a,b,c\in W_9$, $a\cdot b\le c$ iff $b\le a\to c$.
\item  $\bfm{W_9}$ is local, directly indecomposable and $Rad(\bfm{W_9})=\{\top,1,2,3\}$ is its unique proper implicative filter.
\item For any $n\ge 2$, $3.x^n$ is radical term for $\bfm{W_9}$.
\item For any $a\in W_9$,  $a\notin Rad(\bfm{W_9})$ iff $a^3=0$, and so for any $m>2$, $3.x^m$ is Boolean term for $\bfm{W_9}$ and so for $HSP(\bfm{W_9})$.
\item $\bfm{W_9}$ has, up isomorphism, five subalgebras:  it self $\bfm{W_9}$ , the two element Boolean algebra $\bfm{L}_2$ with universe $\{\bot,\top\}$,$\bfm{W_8}$ whose  universe is $W_8=\{\top,1,3,4,5,6,7,\bot\}$,   $\bfm{W_7}$ whose  universe is $W_7=\{\top,1,2,3,6,7\bot\}$ and $\bfm{W_6}$ whose universe is $W_6=\{\top,1,3,6,7\bot\}$.

 If  $8\le k\le 9$, then $\bfm{W_k}$ has two non-trivial homomorphic images: $\bfm{L_3}$ (see above example) and itself; and if  $6\le k\le 7$, then $\bfm{W_k}$ also has two non-trivial homomorphic images: $\bfm{L_2}$  and itself.
 Therefore, since $\bfm{W_9}$ is finite and its gene\-rated variety is congruence distributive by Corollary III 6.10 of \cite{BS81}, we have
 \[HSP(\bfm{W_9})= SP\{\bfm{L}_2,\bfm{L_3},\bfm{W_6},\bfm{W_7},\bfm{W_8},\bfm{W_9}\},\]
 and so, by Lemmata \ref{L: SbProdH} and \ref{L: SbProdH2}, for any $m>2$ $3.x^m$ is also radical term for $HSP(\bfm{W_9})$.
 \item Finally, observe that $HSP(\bfm{W_9})_{\mbb{SS}}=HSP(\bfm{L_3})$.
\ei} \e{example}

\bigskip

\no {\sc Antoni Torrens i Torrell}\\
 Universitat de Barcelona, Spain\\
\texttt{atorrens@ub.edu}

\end{document}